\documentclass[twoside]{article}
\usepackage{eurosym}
\usepackage[utf8]{inputenc}
\usepackage{graphicx}
\usepackage{latexsym}
\usepackage{amsmath}
\usepackage{amsfonts}
\usepackage{amssymb}
\usepackage{hyperref}
\usepackage[usenames]{color}

\newtheorem{theorem}{Theorem}[section]
\newtheorem{corollary}[theorem]{Corollary}
\newtheorem{lemma}[theorem]{Lemma}

\newtheorem{remark}[theorem]{Remark}
\newtheorem{example}[theorem]{Example}

\input{tcilatex}
\begin{document}

\title{\sffamily\bfseries Beyond the classical Strong Maximum Principle:
sign-changing forcing term and flat solutions}
\author{Jes\'{u}s Ildefonso D\'{\i}az and Jes\'{u}s Hern\'{a}ndez\thanks{%
The authors are partially supported by the project PID2020-112517GB-I00 of
the DGISPI, Spain. \hfil\break \indent {\sc
		Keywords}: Strong maximum principle, forcing changing sign, positive flat
solutions, unique continuation, heat equation. \hfil\break 
}}
\date{}
\maketitle

\begin{abstract}
{\small We show that the classical Strong Maximum Principle, concerning
positive supersolutions of linear elliptic equations vanishing on the
boundary of the domain can be extended, under suitable conditions, to the
case in which the forcing term is sign-changing. In addition, for the case
of solutions the normal derivative on the boundary may also vanish on the
boundary (definition of flat solution). This leads to examples in which the
unique continuation property fails. As a first application, we show the
existence of positive solutions for a sublinear semilinear elliptic problem
of indefinite sign. A second application, concerning the positivity of
solutions of the linear heat equation, for some large values of time, with
forcing and/or initial datum changing sign is also given. }
\end{abstract}

\section{Introduction}

In a pioneering paper, on 1910, S. Zaremba \cite{Zaremba1910} established
the well-known \textit{Strong Maximum Principle }saying, in a simple
formulation, that if $\Omega $ is a smooth bounded domain in $\mathbb{R}^{N}$
and a function $u$ verifies 
\begin{equation}
\left\{ 
\begin{array}{ll}
-\Delta u\geq f(x) & \text{in }\Omega , \\ 
u=0 & \text{on }\partial \Omega ,%
\end{array}%
\right.  \label{Supersolution problem}
\end{equation}%
\noindent then 
\begin{equation}
u(x)>0\text{ in }\Omega ,  \tag{P$_{u}$}  \label{Positiv u}
\end{equation}

\noindent assumed that 
\begin{equation}
f(x)\geq 0\text{ in }\Omega ,\text{ }f\neq 0.  \tag{P$_{f}$}
\label{Positiv f}
\end{equation}%
The extension to a more general second order elliptic operator was due to E.
Hopf, on 1927, in his famous paper \cite{Hopf}. Moreover, some years later,
on 1952, E. Hopf \ \cite{Hopf1952} and O.A. Oleinik \cite{Oleinik1952},
independently, proved that under the above conditions the normal derivative
of $u$ satisfies the following sign condition

\begin{equation}
\frac{\partial u}{\partial n}<0\text{ on }\partial \Omega ,  \label{Hypo2.3}
\end{equation}%
(see, the survey \cite{Apushkinskaya-Nazarov} for many other historical
informations).

A \textquotedblleft quantitative strong maximum\textquotedblright\ principle
was obtained since 1987 (\cite{Morel-Ostwald}, \cite{DMO},\cite{Zhao}, \cite%
{Brez-Cabre}, \cite{BerdanDRako},...). In order to be more precise, we will
work in the class of very weak supersolutions: the theory of very weak
solutions of (\ref{ProbLineal}) below was introduced in an unpublished paper
by Ha\"{\i}m Brezis on 1971, later reproduced in (\cite{Brez-Caz et al})
(see a regularity extension in \cite{Diaz-Rako}). This theory applies to the
more general class of data $f\in L_{loc}^{1}(\Omega )$ for which it is
possible to give a meaning to the notion of solution of the corresponding
problem. We assume 
\begin{equation*}
f\in L^{1}(\Omega :\delta )=\left\{ g\in L_{loc}^{1}(\Omega )\text{ such that%
}\int_{\Omega }\left\vert g(x)\right\vert \delta (x)dx<\infty \right\} ,
\end{equation*}%
where 
\begin{equation*}
\delta (x)=d(x,\partial \Omega ).
\end{equation*}%
Then, the above mentioned \textquotedblleft quantitative strong maximum
principle\textquotedblright\ says that if $f$ satisfies (\ref{Positiv f}),
and thus 
\begin{equation}
\int_{\Omega }f(x)\delta (x)dx>0,  \label{strictPositivity}
\end{equation}%
then we get the following estimate, called in \cite{BerdanDRako} as the 
\textit{Uniform Hopf Inequality} (UHI), 
\begin{equation}
u(x)\geq C\delta (x)\int_{\Omega }f(x)\delta (x)dx\text{ a.e. }x\in \Omega ,
\label{distance}
\end{equation}%
for some $C>0$ only dependent on $\Omega $. Notice that, if for instance $%
u\in W^{1,1}(\Omega )$, this implies $\frac{\partial u}{\partial n}<0$ on $%
\partial \Omega .$

The main goal of this paper is to show that the sign assumption (\ref%
{Positiv f}) can be removed so that, under suitable conditions, any
supersolution $u$ satisfying (\ref{Supersolution problem}) for suitable
sign-changing functions $f(x)$ is again strictly positive on $\Omega $.
Moreover, under suitable conditions, this strictly positive supersolution
for some changing sign datum $f(x)$ does not satisfy the condition (\ref%
{Hypo2.3}). In some cases this kind of sign-changing datum $f(x)$ still may
satisfy the condition (\ref{strictPositivity}) (see, Remark \ref{Rem
Integral positive}) but the conclusion (\ref{distance}) may fail (see the
notion of \textit{flat solution} given below).

As far as we know, curiously enough, such type of extension of the classical
strong maximum principle was not presented in the previous literature on the
subject (see, e.g., \cite{Protter-Weimberger}, \cite{Pucci-Serrin}, \cite%
{Vazquez}, \cite{Brezis-Nirenberg} and \cite{Brezis-Ponce}, among many other
papers and books: the survey \cite{Apushkinskaya-Nazarov} contains more than
230 references on the subject, until 2022, but it seems that none of them
deals with the case in which $f(x)$ changes sign). There are some papers in
the literature which could lead to some related conclusion but their
statements are not presented in the same form that in this paper (see, e.g.,
Remark \ref{Rem Siravov} below).

We will pay a special attention to the case in which $f(x)<0$ in some
neighborhood of $\partial \Omega $, but many other cases can be also
considered (see Remark \ref{Rem f negative interior})$.$ In order to present
our results we will use the decomposition%
\begin{equation*}
f(x)=f^{+}(x)-f^{-}(x)
\end{equation*}%
with%
\begin{equation*}
f^{+}(x)=\max (f(x),0)\text{, }f^{-}(x)=-\min (f(x),0)
\end{equation*}%
(notice that $f^{-}(x)\geq 0$). We assume in that paper that the region
where $f(x)$ is negative has at least a part which is touching $\partial
\Omega $, nevertheles different cases can be also considered by similar
arguments (see Remark \ref{Rm First order} below). We assume that there
exists an open subset $\Omega ^{+}\subset \Omega $ such that 
\begin{equation}
\left\{ 
\begin{array}{ll}
f(x)\geq 0 & \text{a.e. }x\in \Omega ^{+}, \\ 
f(x)\leq 0 & \text{a.e. }x\in \Omega \setminus \Omega ^{+}, \\ 
\begin{array}{l}
\partial \Omega \text{ satisfies the interior sphere condition,} \\ 
\sup_{x\in \Omega ^{+}}d(x,\partial \Omega )>0.%
\end{array}
& 
\end{array}%
\right.  \label{hypo fmas omega}
\end{equation}%
The kind of \textquotedblleft new assumptions\textquotedblright\ on $f(x)$
giving positive solutions are of the following type:

\begin{quote}
(H$_{1}$) \textit{A suitable balance expressing that the negative zone of }$%
f(x)$\textit{\ takes place near the boundary }$\partial \Omega $: condition (%
\ref{hypo fmas omega}) holds and there exists a compact set $K\subset \Omega
^{+}$ where $f\neq 0$ on $K$ and 
\begin{equation}
\int_{\Omega }f^{+}\delta >\frac{C_{K}^{\ast }}{c_{^{\ast }K}}\int_{\Omega
}f^{-}\delta ,  \label{Hypo balance f and K}
\end{equation}%
with $c_{K}^{\ast }$ $<C_{K}^{\ast },$ some positive constants (depending on 
$K$) which will be defined later (see Lemma \ref{Lemm uniform y} below).
\end{quote}

\bigskip

\noindent \lbrack This will allow us to conclude that $u\geq C^{+}$ on $%
\partial K,$ for a suitable $C^{+}>0$: see expression (\ref{Positivity
estimate frontera K}) below],

\bigskip

\begin{quote}
\noindent (H$_{2}$) \textit{A suitable decay of }$f(x)$ \textit{near the
boundary }$\partial \Omega $: there exists $\alpha >1$ such that 
\begin{equation}
\left\{ 
\begin{array}{cc}
\min_{\overline{\Omega -K}}((\alpha -1)\left\vert \nabla \varphi
_{1}\right\vert ^{2}-\lambda _{1}\varphi _{1}^{2})>0 & \text{and } \\[0.2cm] 
f(x)\geq -M\varphi _{1}(x)^{\alpha -2}\text{ a.e. }x\in \Omega , & 
\end{array}%
\right.  \label{Hypo f near boundary}
\end{equation}%
with $K$ the compact mentioned in (H$_{1}$) and 
\begin{equation}
M=\frac{\alpha C^{+}\min_{\overline{\Omega -K}}((\alpha -1)\left\vert \nabla
\varphi _{1}\right\vert ^{2}-\lambda _{1}\varphi _{1}^{2})}{\max_{\partial
K}\varphi _{1}{}^{\alpha }}.  \label{Choice M}
\end{equation}
\end{quote}

\noindent Here $\varphi _{1}$ denotes the first eigenfunction of the
Laplacian operator in $\Omega ,$ given by 
\begin{equation}
\left\{ 
\begin{array}{ll}
-\Delta \varphi _{1}=\lambda _{1}\varphi _{1} & \text{in }\Omega , \\ 
\varphi _{1}=0 & \text{on }\partial \Omega ,%
\end{array}%
\right.  \label{First eingenfunction}
\end{equation}%
with $\varphi _{1}>0$ normalized by $\left\Vert \varphi _{1}\right\Vert
_{L^{\infty }(\Omega )}=1$. We recall that by well-known results $\varphi
_{1}\sim \delta $, $\delta (x)=d(x,\partial \Omega )$. Notice that the first
condition in (H$_{2}$) has a geometrical meaning (see Remark~\ref{Rm
Geometrical )}).

Under such type of \textquotedblleft new assumptions\textquotedblright\ on $%
f $ \ we will prove:

\noindent (A) the positivity of $u,$ property (\ref{Positiv u}), still
holds. In addition, if, for instance, $u\in W^{1,1}(\Omega )$ then $\frac{%
\partial u}{\partial n}\leq 0$ on~$\partial \Omega ,$

\noindent (B) under additional conditions on $f(x)$, the positive \textit{%
solution} of the linear problem%
\begin{equation}
\left\{ 
\begin{array}{ll}
-\Delta u=f(x) & \text{in }\Omega , \\ 
u=0 & \text{on }\partial \Omega ,%
\end{array}%
\right.  \label{ProbLineal}
\end{equation}%
[i.e., now with the equality symbol $=,$ instead $\geq $] does not satisfy (%
\ref{Hypo2.3}) but $\frac{\partial u}{\partial n}=0$ on $\partial \Omega $ .

Property (B) corresponds to the notion of \textit{flat solution} already
considered by different authors in the framework of some nonlinear problems
(see, e.g., \cite{D-Pitman}, \cite{Ilyasov-Egorov}, \cite{DHIlyasov2015},
etc.). The existence of flat solutions shows that assumption (\ref{Positiv f}%
) is necessary to conclude (\ref{distance}). Notice also that a flat
solution $u$ on a problem (\ref{ProbLineal}) on the domain $\Omega $ can be
extended by zero to get the unique solution $\widetilde{u}$ of a similar
problem associated to an extended domain $\widetilde{\Omega }\supsetneq
\Omega $ with the right hand side given by 
\begin{equation*}
\widetilde{f}(x)=\left\{ 
\begin{array}{ll}
f(x) & \text{if }x\in \Omega , \\ 
0 & \text{if }x\in \widetilde{\Omega }-\Omega .%
\end{array}%
\right.
\end{equation*}%
In this way we can construct \textit{solutions with compact support} for
data with compact support becoming negative near the boundary of its
support. This proves that the version of the strong maximum principle
obtained in \cite{Brezis-Ponce} (ensuring that the solution $u\geq 0$ of a
linear problem (\ref{ProbLineal}) corresponding to a datum $f\geq 0$, cannot
vanish on some positively measured subset of $\Omega $ except if $u\equiv 0$
on $\Omega $) has optimal conditions on $f(x).$

It is a curious fact that the above considerations are motivated, in some
sense, after the long experience in the study in different semilinear
equations with a non-Lipschitz perturbation in the last fifty years. For
instance, in the study of semilinear problems%
\begin{equation}
\left\{ 
\begin{array}{ll}
-\Delta u=g(x)-\beta (u) & \text{in }\Omega , \\ 
u=0 & \text{on }\partial \Omega ,%
\end{array}%
\right.  \label{PbSemilinear}
\end{equation}%
with $g\geq 0$ and $\beta $ a continuous function, for instance, such that $%
\beta (0)=0$, it is well known the existence of a flat solution under
suitable conditions on $\beta $ and $g$. For the case of $\beta $
non-decreasing, subdifferential of the convex function $j$, $\beta =\partial
j$, such that 
\begin{equation*}
\int_{0}\frac{ds}{\sqrt{j(s)}}<+\infty
\end{equation*}%
and $g\neq 0$ we send the reader to Theorem 1.16 of \cite{D-Pitman} (the so
called \textquotedblleft non-diffusion of the support
property\textquotedblright : see also the generalization presented in \cite%
{brezis ushpeki} and \cite{Alvarez-D} for the case of $\beta $ a multivalued
maximal monotone graph). For the autonomous case $g\equiv 0$ and $\beta $
non monotone see, e.g., \cite{DH-CRAS}, \cite{DHIliasov}, \cite{Diaz Hern
Ilys} and its references. This means that if we take 
\begin{equation}
f(x)=g(x)-\beta (u(x))  \label{Particular f}
\end{equation}%
with $u$ the flat solution of (\ref{PbSemilinear}) then $u$ is also a flat
solution of the corresponding linear problem~(\ref{ProbLineal}). Note that,
necessarily, such $f(x)$ becomes negative near the boundary $\partial \Omega
.$

The above extension of the strong maximum principle admits many
generalizations which will be indicated in form of a series of remarks (Schr%
\"{o}dinger equation [Remark \ref{Rm Schroedinger}], operators with a first
order term [Remark \ref{Rm First order}], linear non-local operators [Remark %
\ref{Rm Nonlocal}], nonlinear elliptic operators and obstacle problem
[Remark \ref{Rem p-Laplace and obstacle problem}], etc.).

The organization of this paper is the following: we start by proving, in
Section 2, a version, as simple as possible, for the one-dimensional case in
which assumptions (H$_{1}$) and (H$_{2}$) \ can be easily formulated in an
optimal way. The $N$-dimensional case and without symmetry conditions is
presented in Section 3. An application to some \textit{sublinear} \textit{%
indefinite} semilinear equations (see, e.g., \cite{Hernandez-Mancebo-Vega})
and \cite{DiazHern-Ilyas2023})) will be given in Section 4. Finally, in
Section 5 we will consider the linear parabolic problem%
\begin{equation*}
\left\{ 
\begin{array}{ll}
u_{t}-\Delta u=f(x,t) & \text{in }\Omega \times (0,+\infty ), \\ 
u=0 & \text{on }\partial \Omega \times (0,+\infty ), \\ 
u(x,0)=u_{0}(x) & \text{on }\Omega .%
\end{array}%
\right.
\end{equation*}%
We will show that the above arguments for stationary equations, jointly to
some related results (\cite{DF}), allow to prove that, under suitably
changing sign conditions on $u_{0}(x)$ and/or on $f(x,t)$ (even if it occurs
for any $t>0$), we get the global positivity of $u(x,t)$ on $\Omega $, for
large values of time ( $\exists $ $t_{0}\geq 0$ such that $u(x,t)>0$ a.e. $%
x\in \Omega ,$ for any $t>t_{0}$).

\section{The symmetric one-dimensional linear problem}

For the sake of the exposition, here we consider supersolutions $u(x)$ of
the symmetric one-dimensional linear problem on the domain $\Omega =(-R,R)$%
\begin{equation}
\left\{ 
\begin{array}{lr}
-u^{\prime \prime }(x)\geq f(x) & \text{in (}-R,R), \\ 
u(\pm R)=0. & 
\end{array}%
\right.  \label{Supersolution unidimens}
\end{equation}%
We assume the symmetry condition 
\begin{equation*}
\text{ }f(x)=f(-x),\text{ }f=f^{+}-f^{-},
\end{equation*}%
and we will work in the framework of the space $L^{1}(\Omega :\delta )$,
with $\delta (x)=d(x,\partial \Omega )$ (and then $\delta (r)=R-r$ if $r\in
(0,R)$)$.$ We assume 
\begin{equation}
f\in L^{1}(\Omega :\delta )\text{, i.e. }\int_{0}^{R}\left\vert
f(s)\right\vert (R-s)ds<\infty ,  \label{Hypo f distancia borde}
\end{equation}%
and we consider very weak supersolutions, i.e., functions such that $u\in
L^{1}(-R,R)$, with $u^{\prime \prime }\in $ $L^{1}(\Omega :\delta )$,
satisfying 
\begin{equation}
-\int_{\Omega }u\psi ^{\prime \prime }\geq \int_{\Omega }f\psi
\label{def very weak}
\end{equation}%
for any $\psi \in W^{2,\infty }(\Omega )\cap W_{0}^{1,\infty }(\Omega )$
such that $\psi \geq 0$. Notice that since any function $\psi \in
W^{2,\infty }(\Omega )\cap W_{0}^{1,\infty }(\Omega )$ satisfies that $%
\left\vert \psi (x)\right\vert \leq C\delta (x)$ for any $x\in \overline{%
\Omega }$, for some $C>0$, then the expressions in (\ref{def very weak})
make sense. The notion of very weak solution is similar but replacing the
symbol $\geq $ by $=$. In some parts of our exposition we will refer to
symmetric solutions (and not merely supersolutions). By well-known results
(see, e.g. \cite{Brez-Caz et al}, \cite{Diaz-Rako}) we have that $u\in
C([0,R])\cap C^{1}[0,R)$, \textit{\ }$u=u(r)$, $r=\left\vert x\right\vert $ 
\begin{equation}
\left\{ 
\begin{array}{lr}
-u^{\prime \prime }(r)=f(r) & \text{in }(0,R), \\ 
u(R)=0,\text{ }u^{\prime }(0)=0. & 
\end{array}%
\right.  \label{Linear onedimens equation}
\end{equation}%
Let us see how the type of assumptions (H$_{1}$) and (H$_{2}$), mentioned in
the Introduction,\ can be easily formulated, and even in an optimal way.

\begin{theorem}
\textbf{\label{Thm onedimensional}} \textit{We} \textit{assume that }$f(x)$
becomes negative near the boundary in the following sense: \textit{there
exists }%
\begin{equation}
r_{0}\in (0,R),
\end{equation}%
\textit{such that}%
\begin{equation}
\left\{ 
\begin{array}{cc}
f(x)=f^{+}(x)\geq 0\text{ if }x\in (0,r_{0}),\text{ }f^{+}\neq 0\text{ on }%
(0,r_{0}), &  \\[0.2cm] 
f(x)=-f^{-}(x)\leq 0\text{ if }x\in (r_{0},R), & 
\end{array}%
\right.  \label{Hypo f+ non zer}
\end{equation}

\noindent (A) Assume the \textquotedblleft balance
condition\textquotedblright 
\begin{equation}
\int_{0}^{r_{0}}f^{+}(s)(R-r_{0})ds>\int_{r_{0}}^{R}f^{-}(s)(R-s)ds
\label{Hypo positivity r0}
\end{equation}%
and the \textquotedblleft decay condition\textquotedblright \textit{\ }%
\begin{equation}
\int_{r}^{R}\left (\int_{r_{0}}^{t}f^{-}(s)ds\right
)dt<(R-r)\int_{0}^{r_{0}}f^{+}(s)ds, \text{ \textit{for any }}r\in (r_{0},R).
\label{Hypo r0 doble}
\end{equation}%
\textit{Then any symmetric supersolution }$u$ satisfies\textit{\ }%
\begin{equation}
u>0\text{ in }(-R,R)\text{\textit{.} }  \label{Charact u positiv}
\end{equation}%
\textit{Moreover, }if for instance $u\in C^{1}[0,R]$, then \textit{we have }%
\begin{equation}
u^{\prime }(R)\leq 0\text{ and }u^{\prime }(-R)\geq 0.
\label{property u' boundary}
\end{equation}

\noindent In addition, assumed (\ref{Hypo positivity r0}), if $u$ is a
solution then $u>0$ \textit{if and only if the decay condition (\ref{Hypo r0
doble}) holds. }

\noindent (B) \textit{Assume now (\ref{Hypo positivity r0}), (\ref{Hypo r0
doble}) and }%
\begin{equation}
f^{-}\in L^{1}(\Omega )\text{, i.e., }\int_{r_{0}}^{R}f^{-}(s)ds<\infty .
\label{Hypo f- L^1}
\end{equation}

\noindent \textit{\ Let }$u$ be \textit{the unique solution }$u$ of (\ref%
{Linear onedimens equation}). Then $u$ \textit{is flat} ($u^{\prime }(\pm
R)=0$) \textit{if and only if the following condition holds }%
\begin{equation}
\int_{0}^{r_{0}}f^{+}(s)ds=\int_{r_{0}}^{R}f^{-}(s)ds\text{, i.e. }
\int_{0}^{R}f(s)ds=0.  \label{Hypo equal integrals}
\end{equation}
\end{theorem}

\begin{remark}
\textrm{\label{Rem Integral positive}It is easy to see that assumptions (\ref%
{Hypo f+ non zer}) and (\ref{Hypo positivity r0}) imply that $f(x)$
satisfies the positivity of the weighted integral (\ref{strictPositivity}).
Indeed, 
\begin{equation*}
\int_{0}^{r_{0}}f(s)(R-s)ds\geq
\int_{0}^{r_{0}}f^{+}(s)(R-r_{0})ds>\int_{r_{0}}^{R}f^{-}(s)(R-s)ds.
\end{equation*}%
\noindent As a matter of fact, one of the main goals of the present paper is
to prove a conjecture raised by the first author and communicated to Jean
Michel Morel in 1985 (when preparing the joint paper \cite{DMO}), concerning
the possibility for sign-changing data $f(x)$ to keep the positivity of the
weighted integral~(\ref{strictPositivity}) and guaranteeing also the
positivity of the supersolution.$\hfill \mbox{$\quad{}_{\Box}$}$ }
\end{remark}

\begin{example}
\textrm{Before to give the proof, let us see the different kinds of
behaviour of solutions for a simple example with a changing sign right hand
side term. Consider 
\begin{equation}
\left\{ 
\begin{array}{lr}
-u^{\prime \prime }(x)=f(x) & \text{in }(-a,a), \\ 
u(\pm a)=0, & 
\end{array}%
\right.  \label{probExam1}
\end{equation}%
for different values of $a\geq 1$ with%
\begin{equation}
f(x)=\left\{ 
\begin{array}{ll}
1 & \text{if }x\in (-1,1), \\ 
-1 & \text{if }x\in (-a,-1)\cup (1,a).%
\end{array}%
\right.  \label{ForcingExamp}
\end{equation}
}
\begin{figure}[h]
	\begin{center}
		\includegraphics[height=13cm,width=16cm]{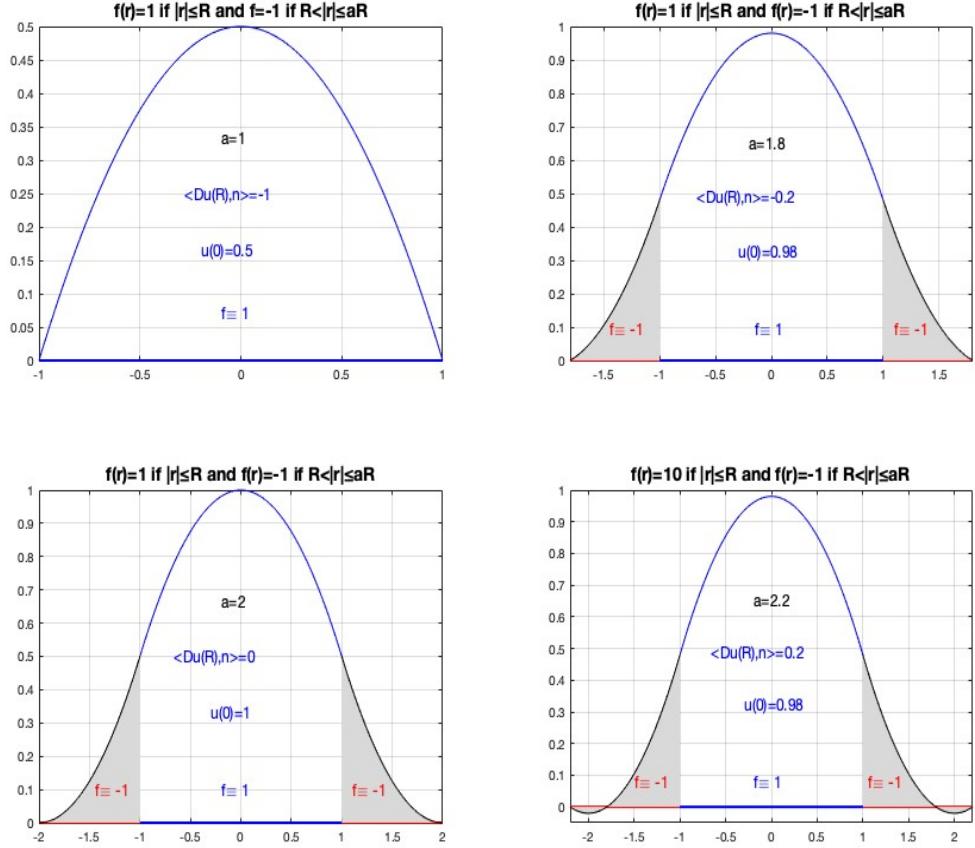}
		\vspace*{-1cm}
		\caption{\small Representation of the exact solution of (\protect\ref{probExam1}),
			and the value of its normal derivative at the boundary, when the forcing is
			given by (\protect\ref{ForcingExamp}), for different values of $a$.}
	\end{center}
\end{figure}
\textrm{\noindent The Figure 1 below shows different cases for the values of 
$a=1$, \ $a=1.8$, $a=2$ and $a=2.2.$ }

\textrm{\noindent In the first case of Figure 1 $a=1$, the forcing $f(x)$ is
strictly positive and the classical strong maximum principle applies. In the
case $a=1.8$ we see that the solution is strictly positive and its normal
derivative at the boundary is again strictly negative. For $a=2$ the
conditions of part B of Theorem \ref{Thm onedimensional} hold and we see
that the solution is flat. Finally, for $a=2.2$ the assumptions of part A of
Theorem \ref{Thm onedimensional} fail and the solution becomes negative near
the boundary. $\hfill\mbox{$\quad{}_{\Box}$}$ }
\end{example}

\bigskip \noindent \textit{Proof of Theorem \ref{Thm onedimensional}.} By
well-known results (see, e.g., \cite{Brez-Caz et al}, \cite{Brez-Cabre}, 
\cite{Diaz-Rako}) we may assume that\textit{\ }$u$ is\ symmetric\ and has
some regularity properties, i.e. \ $u$ is such that $u\in C([0,R])\cap
C^{1}([0,R))$ with $u^{\prime \prime }\in L^{1}((0,R):R-r)$, \textit{\ }$%
u=u(r)$, $r=\left\vert x\right\vert $ and%
\begin{equation}
\left\{ 
\begin{array}{lr}
-u^{\prime \prime }(r)\geq f(r) & \text{in (}0,R), \\ 
u(R)=0,u^{\prime }(0)=0. & 
\end{array}%
\right.
\end{equation}%
To prove part (A) let us start by proving that 
\begin{equation}
\text{(\ref{Hypo positivity r0}) implies that }u(r_{0})>0\text{.}
\label{Chracter u(r_0)}
\end{equation}%
Multiplying the equation by $(R-r)$, with $r\in (r_{0},R),$ and using (\ref%
{Hypo f+ non zer}) we have 
\begin{equation*}
-\int_{r}^{R}u^{\prime \prime }(s)(R-s)ds\geq -\int_{r}^{R}f^{-}(s)(R-s)ds.
\end{equation*}%
Integrating by parts, and using the boundary condition we get%
\begin{equation*}
-u^{\prime }(R)(R-R)+u^{\prime }(r)(R-r)-\int_{r}^{R}u^{\prime }(s)ds\geq
-\int_{r}^{R}f^{-}(s)(R-s)ds,
\end{equation*}%
i.e., since $u(R)=0,$ for any $r\in (r_{0},R)$ we get 
\begin{equation}
u^{\prime }(r)(R-r)+u(r)\geq -\int_{r}^{R}f^{-}(s)(R-s)ds.
\label{Balance near boundary}
\end{equation}%
On the other hand, from (\ref{Hypo f distancia borde}) and the structure
assumptions (\ref{Hypo f+ non zer}) we have that in fact $f^{+}\in
L^{1}(0,r_{0}).$ Then, for any $r\in \lbrack 0,r_{0}]$ we have%
\begin{equation}
\int_{0}^{r}u^{\prime \prime }(s)ds=u^{\prime }(r)\leq
-\int_{0}^{r}f^{+}(s)ds\leq 0  \label{Integral near 0}
\end{equation}%
(which clearly implies that the maximum of $u$ is taken at $r=0$). Thus,
from (\ref{Hypo f+ non zer}) 
\begin{equation}
u^{\prime }(r_{0})\leq -\int_{0}^{r_{0}}f^{+}(s)ds<0.
\label{property u'(r_0)}
\end{equation}%
Substituting in (\ref{Balance near boundary}), since $u\in C([0,R])\cap
C^{1}(0,R),$ we get that 
\begin{equation*}
-(R-r_{0})\int_{0}^{r_{0}}f^{+}(s)ds+u(r_{0})\geq
-\int_{r_{0}}^{R}f^{-}(s)(R-s)ds,
\end{equation*}%
which proves that (\ref{Hypo positivity r0}) implies property (\ref{Chracter
u(r_0)}), i.e., $u(r_{0})>0$. [In fact, it is easy to see that if $u$ is a
solution (and not merely a supersolution) then assumption (\ref{Hypo
positivity r0}) is also a necessary condition to have $u(r_{0})>0$]. Notice
also that, from (\ref{Hypo f+ non zer}) we always have%
\begin{equation*}
\int_{0}^{r_{0}}f^{+}(s)(R-s)ds>(R-r_{0})\int_{0}^{r_{0}}f^{+}(s)ds.
\end{equation*}

\noindent Then we get 
\begin{equation}
\left\{ 
\begin{array}{lr}
u^{\prime \prime }(r)\leq -f^{+}(r)\leq 0 & \text{in (}0,r_{0}), \\ 
u(r_{0})>0,u^{\prime }(0)=0, & 
\end{array}%
\right.
\end{equation}%
which implies, from (\ref{Integral near 0}), that 
\begin{equation*}
u(r)\geq u(r_{0})>0\text{ for any }r\in \lbrack 0,r_{0}].
\end{equation*}%
\noindent To complete the proof of (A), we see that from the structure
conditions (\ref{Hypo f+ non zer}), for any $r\in \lbrack r_{0},R)$ we have 
\begin{equation}
\int_{r_{0}}^{r}u^{\prime \prime }(s)ds\leq \int_{r_{0}}^{r}f^{-}(s)ds.
\label{property u' nonpositive}
\end{equation}%
Then, integrating again, for any $\widehat{r}\in (r,R)$ we have

\begin{equation*}
-\int_{r}^{\widehat{r}}u^{\prime }(t)dt+u^{\prime }(r_{0})(\widehat{r}%
-r)\geq -\int_{r}^{\widehat{r}}\left (\int_{r_{0}}^{t}f^{-}(s)ds\right )dt,
\end{equation*}%
i.e., from (\ref{property u'(r_0)})%
\begin{equation}
-u(\widehat{r})+u(r)-(\widehat{r}-r)\int_{0}^{r_{0}}f^{+}(s)ds\geq
-\int_{r}^{\widehat{r}}\left (\int_{r_{0}}^{t}f^{-}(s)ds\right )dt.
\label{Balance r^}
\end{equation}%
Assumption (\ref{Hypo r0 doble}) implies that, 
\begin{equation}
\int_{0}^{R}(\int_{r_{0}}^{r}f^{-}(s)ds)dr<+\infty ,
\label{Proeperty double integral}
\end{equation}%
and thus, for any $r\in \lbrack r_{0},R)$ 
\begin{equation*}
\underset{\widehat{r}\nearrow R}{\lim }\int_{r}^{\widehat{r}}\left
(\int_{r_{0}}^{t}f^{-}(s)ds\right
)dt=\int_{r}^{R}(\int_{r_{0}}^{t}f^{-}(s)ds)dt<+\infty \text{.}
\end{equation*}%
Then, making $\widehat{r}\nearrow R$, in (\ref{Balance r^}) we get 
\begin{equation}
u(r)\geq (R-r)\int_{0}^{r_{0}}f^{+}(s)ds-\int_{r}^{R}\left
(\int_{r_{0}}^{t}f^{-}(s)ds\right )dt,  \label{Formula explicit}
\end{equation}%
which leads to the positivity conclusion (\ref{Charact u positiv}).

\noindent Obviously, if for instance $u\in C^{1}[0,R]$, since $u(R)=0$ we
conclude that $u^{\prime }(\pm R)\leq 0.$

\noindent To prove (B), we observe that assumed (\ref{Hypo positivity r0})
and (\ref{Hypo r0 doble}), then if $u$ is a solution of (\ref{Linear
onedimens equation}) the inequality~(\ref{property u' nonpositive}) becomes
an equality. Making $r=R$ and using (\ref{property u'(r_0)}), since $u\in
W^{2,1}(0,R)$ we get that $u\in C^{1}[0,R]$ and

\begin{equation*}
u^{\prime }(R)=\int_{r_{0}}^{R}f^{-}(s)ds-\int_{0}^{r_{0}}f^{+}(s)ds.
\end{equation*}%
This proves the necessary and sufficient condition in order to have flat
solutions. $_\Box$

\medskip The existence of nonnegative solutions with compact support in a
larger domain $\widetilde{\Omega }=(-\widetilde{R},\widetilde{R}$), with $%
\widetilde{R}>R$ is a simple consequence of the existence of flat solutions
on the small domain $\Omega =(-R,R)$. Notice that this proves a failure of 
\textit{the unique continuation property} under the below conditions.

\begin{corollary}
\label{Corol compact support onedim}Let $\widetilde{f}\in L^{1}(\widetilde{
\Omega })$ be the extension of a given function $f\in L^{1}(\Omega )$, i.e.
such that 
\begin{equation*}
\widetilde{f}(x)=\left\{ 
\begin{array}{ll}
f(x) & \text{if }x\in \Omega , \\ 
0 & \text{if }x\in \widetilde{\Omega }\setminus\Omega .%
\end{array}
\right.
\end{equation*}%
Assume that $f$ satisfies the conditions (\ref{Hypo positivity r0}), \textit{%
\ (\ref{Hypo r0 doble}) and (\ref{Hypo equal integrals}). Let }$u$ be the
unique solution of (\ref{Linear onedimens equation}) and let $\widetilde{u}$
be the extension of $u$ defined as%
\begin{equation*}
\widetilde{u}(x)=\left\{ 
\begin{array}{ll}
u(x) & \text{if }x\in \Omega , \\ 
0 & \text{if }x\in \widetilde{\Omega }\setminus\Omega .%
\end{array}
\right.
\end{equation*}%
Then $\widetilde{u}\in C([0,\widetilde{R}])\cap C^{1}(0,\widetilde{R})$ and
it is the unique weak solution of the problem 
\begin{equation}
\left\{ 
\begin{array}{lr}
-\widetilde{u}^{\prime \prime }(x)=\widetilde{f}(x) & \text{in }(-\widetilde{
R},\widetilde{R}), \\ 
\widetilde{u}(\pm \widetilde{R})=0. & 
\end{array}
\right.  \label{Linear extended problem}
\end{equation}
\end{corollary}

\noindent \emph{Proof.} By Theorem \ref{Thm onedimensional} we know that $%
u(x)$ is a flat solution of the problem (\ref{Linear onedimens equation})
associated to $f\in L^{1}(\Omega )$ on $\Omega $ and then the extension $%
\widetilde{u}$ is a weak solution of the extended problem (\ref{Linear
extended problem}). By the uniqueness of solutions for such problem, $%
\widetilde{u}(x)$ is the unique function satisfying (\ref{Linear extended
problem}).$\hfill\mbox{$\quad{}_{\Box}$}$

\begin{remark}
\textrm{No flat solution may satisfy the Hopf conclusion (\ref{Hypo2.3}) nor
the decay estimate (\ref{distance}). This proves that the non-negative
condition assumed on $f(x)$ in the corresponding results is necessary.
Analogously, the Corollary \ref{Corol compact support onedim} proves that
the version of the strong maximum principle obtained in \cite{Brezis-Ponce}
(ensuring that the solution $u\geq 0$ of a linear problem (\ref{ProbLineal})
corresponding to a datum $f\geq 0$, cannot vanish on some positively
measured subset of $\Omega $ except if $u\equiv 0$ on $\Omega $) fails for
positive solutions corresponding to changing sign data $f(x).\hfill%
\mbox{$\quad{}_{\Box}$}$ }
\end{remark}

\begin{remark}
\textrm{\label{Rm Rakotoson}It is possible to give an alternative proof of
Theorem \ref{Thm onedimensional} by using the Green function associated to
the Dirichlet problem when the datum $f^{-}(x)$ is symmetric. Indeed, as
before we can argue only for very weak solutions of the problem (\ref%
{probExam1}). Let us assume $a=1.$ Then the Green function $G(x,y)$ is given
by 
\begin{equation}
G(x,y)=\left\{ 
\begin{array}{cc}
\frac{1}{2}(1-y)(1+x) & \text{if }-1<x<y, \\[.2cm] 
\frac{1}{2}(1+y)(1-x) & \text{if }y\leq x<1,%
\end{array}
\right.  \label{Green funct}
\end{equation}
(see, e.g., \cite{Stakgold}, p.54) and since $f\in L^{1}((-1,+1):\delta )$
we know that 
\begin{equation}
u(x)=\int_{-1}^{+1}G(x,y)f(y)dy, \quad x\in \lbrack -1,+1].
\label{Solut by Green}
\end{equation}%
Recall that in this special case $\delta (x)=1-\left\vert x\right\vert $.
Then, some straightforward computations lead to the explicit formula%
\begin{equation}
u(x)=(1-\left\vert x\right\vert )\int_{0}^{\left\vert x\right\vert
}f(s)ds+\int_{\left\vert x\right\vert }^{1}f(s)(1-s)ds,\quad x\in \lbrack
-1,+1]. \text{ }  \label{Explicit radial solut}
\end{equation}%
From this formula, by assuming (\ref{Hypo f+ non zer}), (\ref{Hypo
positivity r0}) and (\ref{Hypo r0 doble}), we can deduce again the
expressions (\ref{Formula explicit}) and~(\ref{Balance r^}), and the proof
follows. Once again, we see that $u>0$ if and only if (\ref{Hypo f+ non zer}%
), (\ref{Hypo positivity r0}) and (\ref{Hypo r0 doble}) hold. $\hfill%
\mbox{$\quad{}_{\Box}$}$ }
\end{remark}

\begin{remark}
\textrm{Clearly, condition (\ref{Hypo f distancia borde}) is weaker than (%
\ref{Hypo f- L^1}). On the other hand, we point out that the decay
assumption (\ref{Hypo r0 doble}) indicates that the indefinite integral $%
\int_{r}^{R}(\int_{r_{0}}^{t}f^{-}(s)ds)dt$ should decay to zero, as $%
r\nearrow R$, with a rate less than a linear decay (and in fact on the whole
interval $(r_{0},R)$)$.$ Notice that this does not require a pointwise decay
of the type $f(r)\rightarrow 0$ as $r\nearrow R$. Indeed, if, for instance, $%
f^{-}(r)=C(R-r)^{-\alpha }$ for some $C>0$ and $\alpha >0$ then condition (%
\ref{Hypo f distancia borde}) holds if and only if $\alpha \in (0,2)$.
Moreover 
\begin{equation*}
\int_{r_{0}}^{t}f^{-}(s)ds=\frac{C}{(1-\alpha )}\big [(R-r_{0})^{1-\alpha
}-(R-t)^{1-\alpha }\big ]
\end{equation*}%
and 
\begin{equation*}
\begin{array}{c}
\displaystyle\int_{r}^{R}\left( \int_{r_{0}}^{t}f^{-}(s)ds\right) dt=\frac{%
C(R-r_{0})^{1-\alpha }(R-r)}{(1-\alpha )}-\frac{C}{(1-\alpha )}%
\int_{r}^{R}(R-t)^{1-\alpha }dt \\[0.2cm] 
\displaystyle=\frac{C(R-r_{0})^{1-\alpha }(R-r)}{(1-\alpha )}+\frac{C}{%
(1-\alpha )(2-\alpha )}(R-r)^{2-\alpha }.%
\end{array}%
\end{equation*}%
Then, the decay condition (\ref{Hypo r0 doble}) holds if we assume $C$ small
and $\alpha \in (0,1)$. Notice that it fails for $\alpha \in (1,2)$ (similar
computations can be carried out for $\alpha =1$).$\hfill \mbox{$\quad{}_{%
\Box}$}$ }
\end{remark}

\begin{remark}
\textrm{Concerning the optimality of the decay condition, as indicated in
Theorem \ref{Thm onedimensional}, in the class of solutions of (\ref{Linear
onedimens equation}), the condition is optimal. It is not difficult to build
explicit examples of functions $f(r)$ satisfying (\ref{Hypo f distancia
borde}) and such that $f(r)>0$ on a very large zone on $r\in (0,r_{0})$ and
only negative in a very small region $r\in (r_{0},R)$, for which the unique
solution $u$ of (\ref{Linear onedimens equation}) is negative near the
boundary $r=R.$ This is the case, for instance, if we consider 
\begin{equation*}
f(r)=\left\{ 
\begin{array}{ll}
F & \text{if }r\in (0,r_{0}) \\ 
-C_{f}(R-r)^{-3/2} & \text{if }r\in (r_{0},R),%
\end{array}%
\right.
\end{equation*}%
for some positive constants $F$ and $C_{f}$. It is clear that $f^{-}\in
L^{1}(\Omega :\delta )$ and thus there exists a unique solution $u$ of the
corresponding problem (\ref{Linear onedimens equation}). Moreover, by using
formula (\ref{Formula explicit}) (which becomes an equality for the case of
solutions, and not merely supersolutions) it is easy to check that, even if $%
\varepsilon =R-r_{0}$ is small, if $F/C_{f}$ is large enough then there
exists $r^{\ast }\in (r_{0},R)$ such that $u(r)>0$ if $r\in \lbrack
0,r^{\ast }),$ $u(r^{\ast })=0$, $u(r)<0$ if $r\in (r^{\ast },R)$ and $%
u(R)=0.\hfill \mbox{$\quad{}_{\Box}$}$ }
\end{remark}

\begin{remark}
\textrm{Notice that, under the structure condition (\ref{Hypo f+ non zer})
on $f(x)$ we have%
\begin{equation*}
\int_{0}^{r_{0}}f^{+}(s)ds=\int_{\Omega }f^{+}(x)dx\text{, }
\int_{r_{0}}^{R}f^{-}(s)ds=\int_{\Omega }f^{-}(x)dx,
\end{equation*}%
and thus condition (\ref{Hypo equal integrals}) is equivalent to 
\begin{equation*}
\int_{0}^{R}f(s)ds=0.
\end{equation*}
$\hfill\mbox{$\quad{}_{\Box}$}$ }
\end{remark}

\begin{remark}
\textrm{\label{Rem f negative interior}The case in which $f>0$ on a large
part of the domain but with $f<0$ in an interior subset of $\Omega $ can be
also considered by this type of techniques. In that case the positive
solutions satisfy that $u^{\prime }(R)<0$ and $u^{\prime }(-R)>0$ but they
may generate interior points $x_{0}\in \Omega $ where $u(x_{0})=u^{\prime
}(x_{0})=0$ and also the formation of an internal \textquotedblleft dead
core\textquotedblright . Notice that, again, this proves a failure of 
\textit{the unique continuation property under such conditions.} Here is an
example with a similar structure to Example 1 but now reversing the positive
and negative subsets of $f(x).$ }

\begin{figure}[h]
	\begin{center}
		\includegraphics[height=12cm,width=16cm]{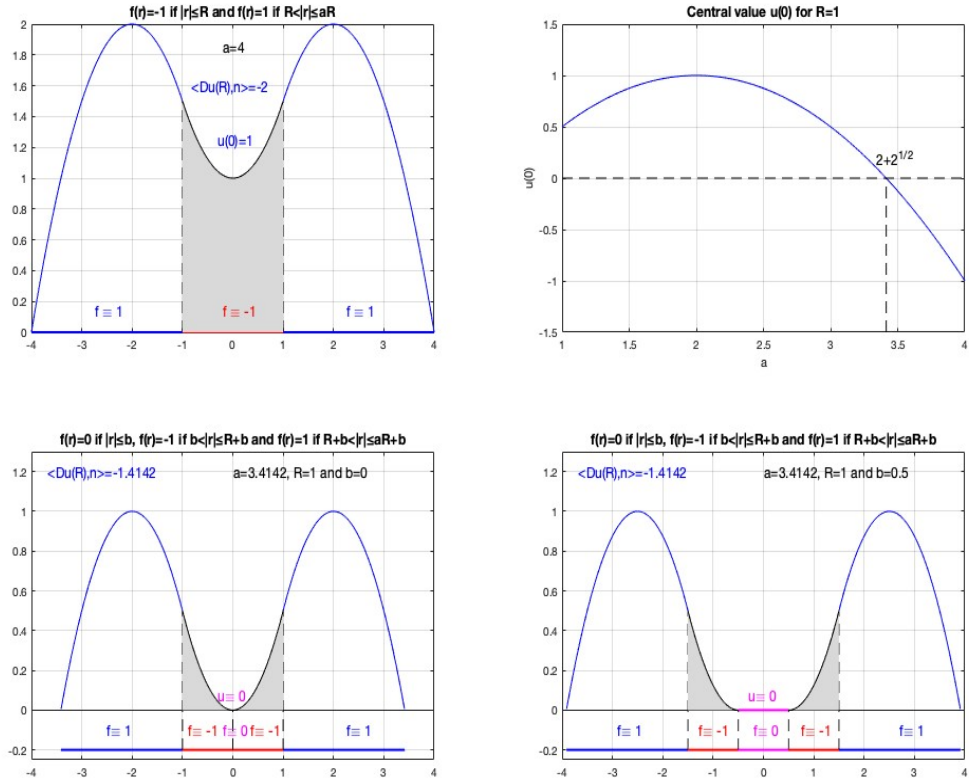}
		\vspace*{-1cm}
		\caption{\small Representation of the solution of problem \protect\ref{probExam1}
			when the forcing is like (\protect\ref{ForcingExamp}) but with the oposite
			sign on $f(x)$.}
	\end{center}
\end{figure}
\textrm{\noindent In the first case of Figure 2 $a=4$, the forcing $f(x)$ is
now positive close to the boundary and negative on $(-1,1)$: the solution is
positive and a local minimum is created at $x=0$. The dependence on $a$ of
the local minimum $u(0)$ is also included in Figure 2. For $a=3.4142$ we
have $u(0)=0$. Finally, In Figure 2 it is also represented the solution $u$
corresponding to a forcing term $f(x)$ vanishing on an interval $(-b,b)$ of $%
x=0$, i.e., 
\begin{equation}
f(x)=\left\{ 
\begin{array}{cl}
1 & \text{if }x\in (-a-b,-1-b)\cup (1+b,a+b). \\ 
0 & \text{if }x\in (-b,b), \\ 
-1 & \text{if }x\in (-1-b,-b)\cup (b,1+b).%
\end{array}
\right.
\end{equation}%
We see the formation of an internal \textquotedblleft dead
core\textquotedblright\ (the interval $(-0.5,0.5)$) where $u=0.\hfill%
\mbox{$\quad{}_{\Box}$}$ }
\end{remark}

\begin{remark}
\textrm{In fact, the property $\frac{\partial u}{\partial n}=0$ is a local
property and it may occur only on a part of $\partial \Omega $. A mixed
situation, exhibiting this fact, is given in the following example (see also 
\cite{Cirmi-Diaz2}). Consider the problem 
\begin{equation}
\left\{ 
\begin{array}{lr}
-u^{\prime \prime }(x)=f(x) & \text{in }(0,1), \\ 
u(0)=u(1)=0, & 
\end{array}
\right.  \label{Example onedim}
\end{equation}%
with 
\begin{equation*}
f(x)=ax-1.
\end{equation*}%
Notice that $f(x)$ changes sign if $a>1.$ It is easy to see that the unique
solution of (\ref{Example onedim}) is%
\begin{equation*}
u(x)=\frac{x^{2}}{2}-\frac{ax^{3}}{6}+\frac{(a-3)x}{6}.
\end{equation*}%
On the other hand, since the roots of $u(x)=0$ are $x_{-}=1$ and $x_{+}=-1+%
\frac{3}{a}$ we see that $u(x)>0$ if $a\geq 3$. In fact, for $a=3$ we have $%
u^{\prime }(0)=0$ and $u^{\prime }(0)>0$ if $a>3.\hfill\mbox{$\quad{}_{%
\Box}$}$ }
\end{remark}

\begin{remark}
\textrm{Notice that a necessary condition in order to have the positivity of
the solutions of \ problem (\ref{Example onedim}) is that 
\begin{equation}
\int_{0}^{1}f(s)\varphi _{1}(s)ds>0,  \label{case A}
\end{equation}%
(it suffices to multiply the equation by $\varphi _{1}$ and integrate twice
by parts). The above example shows that this necessary condition is not
sufficient in order to have the positivity of the solution. Indeed, we can
study the values of $a$ for which $\int_{\Omega }f(x)\varphi _{1}>0$: we have%
\begin{equation*}
\int_{\Omega }f(x)\varphi _{1}(x)dx=\int_{0}^{1}(ax-1)\sin \pi xdx=\frac{a-2 
}{\pi },
\end{equation*}%
and then $\int_{\Omega }f\varphi _{1}>0$ for $a>2$ $\ $and, otherwise, the
positivity of $u$ requires $a\geq 3$ (i.e., for $a\in (2,3]$ we have that $%
\int_{\Omega }f(x)\varphi _{1}>0$ but $u\ngtr 0$). Notice also that if we
take $\Omega =(-1,1),$ 
\begin{equation*}
f(x)=\left\{ 
\begin{array}{ll}
3x-1 & \text{if }x\in (0,1) \\ 
-3x-1 & \text{if }x\in (-1,0)%
\end{array}
\right.
\end{equation*}%
then the solution of 
\begin{equation}
\left\{ 
\begin{array}{lr}
-u^{\prime \prime }(x)=f(x) & \text{in }(-1,1), \\ 
u(-1)=u(1)=0, & 
\end{array}
\right.
\end{equation}%
is%
\begin{equation*}
u(x)=\left\{ 
\begin{array}{ll}
\frac{1}{2}x^{2}(1-x) & \text{if }x\in (0,1) \\[.15cm] 
\frac{1}{2}x^{2}(1+x) & \text{if }x\in (-1,0)%
\end{array}
\right.
\end{equation*}%
and thus, $u(x)>0$ for any $x\in (-1,0)\cup (0,1)$ and $u(0)=0.$ In
addition, if we take 
\begin{equation*}
f(x)=\left\{ 
\begin{array}{ll}
3(x-b)-1 & \text{if }x\in (b,1+b) \\ 
0 & \text{if }x\in (-b,b) \\ 
-3(x+b)-1 & \text{if }x\in (-1-b,-b)%
\end{array}
\right.
\end{equation*}%
for some $b>0$ then the solution has a dead core in $(-b,b),$%
\begin{equation*}
u(x)=\left\{ 
\begin{array}{ll}
\frac{1}{2}(x-b)^{2}(1+b-x) & \text{if }x\in (b,1+b) \\[.15cm] 
0 & \text{if }x\in (-b,b) \\[.15cm] 
\frac{1}{2}(x+b)^{2}(1+b+x) & \text{if }x\in (-1-b,-b).%
\end{array}
\right.
\end{equation*}
$\hfill\mbox{$\quad{}_{\Box}$}$ }
\end{remark}

\section{The N-dimensional case and a general bounded domain $\Omega $}

Now we consider the $N$-dimensional case and a general bounded regular
domain $\Omega $. Our strategy will be quite closed to the main idea of the
proof of Theorem \ref{Thm onedimensional} for the one-dimensional case. We
will assume $f(x)$ such that there exists an open subset $\Omega ^{+}\subset
\Omega $ verifying the hypothesis (\ref{hypo fmas omega}).

In a first step, we will prove something weaker that in the first step of
the proof of Theorem \ref{Thm onedimensional}, we will not prove the
positivity of $u$ on $\Omega ^{+}$ (in Theorem \ref{Thm onedimensional} the
interval $[0,r_{0}]$) but on a regular compact $K$ contained in $\Omega ^{+}$%
: we will show that there exists a positive constant $C^{+}$ such that any
supersolution $u$ of (\ref{Supersolution problem}) satisfies 
\begin{equation}
u\geq C^{+}\text{ on }K,  \label{First step}
\end{equation}%
see the expression (\ref{Positivity estimate frontera K}) below.

\bigskip

For the proof, we will need to work with some auxiliary problems of the type 
\begin{equation}
\left\{ 
\begin{array}{ll}
-\Delta \varsigma _{y}=\chi _{B_{\varrho }(y)} & \text{in }\Omega , \\ 
\varsigma _{y}=0 & \text{on }\partial \Omega ,%
\end{array}
\right.  \label{Auxiliary problems}
\end{equation}%
where $y\in \partial K$, $\chi _{A}$ denotes the characteristic function of $%
A$ and $\varrho >0$ is small enough such that

\begin{equation*}
B_{\varrho }(y)\subset \Omega ^{+}.
\end{equation*}%
Due to the compactness of $K$ and the classical strong maximum principle
(Hopf-Oleinik boundary lemma, since $\partial \Omega $ satisfies the
interior sphere condition), it is well-known that for any $y\in \partial K$
there exist two positive constants $c_{K}(y)<C_{K}(y)$ such that%
\begin{equation}
c_{K}(y)\delta (x)\leq \varsigma _{y}(x)\leq C_{K}(y)\delta (x)\text{ a.e. }
x\in \Omega ,  \label{Inequality distance boundary}
\end{equation}%
where $\delta (x)=d(x,\partial \Omega )$. Using the compactness of $K,$ it
is possible to get the following result (that will be proved later) giving
some uniform estimates:

\begin{lemma}
\label{Lemm uniform y} There exists two positive constants $0<c_{K}^{\ast
}<C_{K}^{\ast }$ such that 
\begin{equation}
c_{K}^{\ast }\delta (x)\leq \varsigma _{y}(x)\leq C_{K}^{\ast }\delta (x) 
\text{ a.e. }x\in \Omega ,\text{ for any }y\in \partial K.
\label{Uniform estimates test functions}
\end{equation}
\end{lemma}

\bigskip

As in the one-dimensional case, there will be a second step in the proof of
the main result of this section, where we prove that, under the balance and
decay conditions mentioned in the Introduction, the unique solution $v$ of
the problem on the ring $\Omega \setminus K$%
\begin{equation}
\left\{ 
\begin{array}{ll}
-\Delta v=-f^{-}(x) & \text{in }\Omega \setminus K, \\ 
v=0 & \text{on }\partial \Omega , \\ 
v=C^{+} & \text{on }\partial K,%
\end{array}
\right.  \label{Problem ring}
\end{equation}%
is a positive subsolution and thus 
\begin{equation*}
0<v(x)\leq u(x)\text{ a.e. }x\in \Omega \setminus K.
\end{equation*}%
For simplicity in the exposition (since there are other different options)
this subsolution $v(x)$ will be constructed in terms of a suitable power of $%
\varphi _{1}$, the normalized first eigenfunction of the Laplacian operator
on $\Omega $. We recall that $u$ is a very weak supersolution and that this
means that $u\in L^{1}(\Omega ),$ and 
\begin{equation}
-\int_{\Omega }u\Delta \psi \geq \int_{\Omega }f\psi
\label{def very weak Ndim}
\end{equation}%
for any $\psi \in W^{2,\infty }(\Omega )\cap W_{0}^{1,\infty }(\Omega )$
such that $\psi \geq 0$. Again, since any function $\psi \in W^{2,\infty
}(\Omega )\cap W_{0}^{1,\infty }(\Omega )$ satisfies that $\left\vert \psi
(x)\right\vert \leq C\delta (x)$ for any $x\in \overline{\Omega }$, for some 
$C>0$, then the expressions in (\ref{def very weak Ndim}) make sense. The
notion of very weak solution is similar but replacing the symbol $\geq $ by $%
=$. We have

\begin{theorem}
\label{Thm N-dimensional}Let $f\in L^{1}(\Omega :\delta ).$ Assume (H$_{1}$)
and (H$_{2}$). Then any supersolution $u$ of (\ref{Supersolution problem})
satisfies that 
\begin{equation*}
u(x)>0\text{ a.e. }x\in \Omega .
\end{equation*}

\noindent If in addition $f\in L^{1}(\Omega )$ and 
\begin{equation}
\int_{\Omega }f(x)dx=0,  \label{Flat solu assumpt}
\end{equation}%
then the unique weak solution $u\in W_{0}^{1,1}(\Omega )$ of the linear
problem (\ref{ProbLineal}) is a flat solution.
\end{theorem}

\noindent \underline{Proof. First step\textit{.}}\textit{\ }As in the proof
of the first part of Lemma 3.2 of \cite{Brez-Cabre} (in which the authors
offer an alternative proof to the main result of \cite{Morel-Ostwald}) we
will use the mean value theorem (in our case on $\Omega ^{+}$) but in a
different way. As mentioned before, let $K$ be a regular compact contained
in $\Omega ^{+}$. Since $K$ is compact and $\Omega ^{+}$ is an open subset
of $\Omega $, $dist(K,\partial \Omega ^{+})>0.$ Let $\varrho
<dist(K,\partial \Omega ^{+}).$Then, 
\begin{equation*}
B_{\varrho }(y)\subset \Omega ^{+}\text{ for any }y\in \partial K\text{, and 
}d(\partial B_{\rho }(y),\partial \Omega ^{+})>0.
\end{equation*}%
Let $\varsigma _{y}(x)$ be the solution of \ the auxiliary problem (\ref%
{Auxiliary problems}). Let us assume, for the moment, that $u\in
C^{0}(\Omega ^{+})$. Then, since $-\Delta u\geq 0$ in $\Omega ^{+}$ we
conclude that there exists a positive constant $\widehat{c}$ such that, for
any $y\in \partial K$ 
\begin{equation*}
u(y)\geq \frac{1}{\left\vert B_{\varrho }(y)\right\vert }\int_{B_{\varrho
}(y)}u=\widehat{c}\int_{\Omega }u(-\Delta \varsigma _{y})
\end{equation*}

\noindent Using (\ref{def very weak Ndim}), i.e. integrating twice by parts
(since $\varsigma _{y}\in W^{2,\infty }(\Omega )\cap W_{0}^{1,\infty
}(\Omega )$ and $\varsigma _{y}$ $\geq 0$), and using the uniform estimates (%
\ref{Uniform estimates test functions}) we have 
\begin{equation*}
\int_{\Omega }u(-\Delta \varsigma _{y})\geq \int_{\Omega }f\varsigma
_{y}\geq c_{K}^{\ast }\int_{\Omega }f^{+}\delta -C_{K}^{\ast }\int_{\Omega
}f^{-}\delta .
\end{equation*}%
Then, thanks to the balance condition (H$_{1}$), we get that 
\begin{equation*}
u\geq C^{+}\text{ on }\partial K,
\end{equation*}%
with%
\begin{equation}
C^{+}=\widehat{c}\left[ c_{K}^{\ast }\int_{\Omega }f^{+}\delta -C_{K}^{\ast
}\int_{\Omega }f^{-}\delta \right] >0.
\label{Positivity estimate frontera K}
\end{equation}%
Moreover, since 
\begin{equation*}
\bigg \{%
\begin{array}{cc}
-\Delta u\geq f\geq 0 & \overset{}{\text{in }\overset{\circ }{K}}, \\%
[0.125cm] 
u\geq C^{+} & \text{on }\partial K,%
\end{array}%
\end{equation*}%
we get (\ref{First step}) when $u\in C^{0}(\Omega ^{+})$.

\noindent In the general case it is enough to approximate $f$ and $\Delta u $
by a sequence of regular functions such that $f_{n}\in L^{\infty }(\Omega ),$%
, $u_{n}\in C^{0}(\Omega )$, satisfying the assumptions of the statement for
each $n\in \mathbb{N}$. We know that the corresponding sequences of
functions $f_{n}$ and $u_{n}$ converge in $L^{1}(\Omega :\delta )$ and $%
L^{1}(\Omega ),$ respectively (see \cite{Brez-Caz et al}, \cite{Diaz-Rako}).
Then we arrive to the conclusion since the inequality (\ref{First step}) is
stable for the strong convergence in $L^{1}(\Omega )$. \bigskip

\noindent \underline{Second step\textit{.}} Let us prove that, under (H$_{1}$%
) and (H$_{2}$), the function 
\begin{equation*}
w(x)=k\varphi _{1}(x)^{\alpha }
\end{equation*}%
is a subsolution to problem (\ref{Problem ring}), for some positive constant 
$k$ and for some $\alpha >1$, where $\varphi _{1}$ denotes the first
eigenfunction of the Laplacian operator in $\Omega $ (see (\ref{First
eingenfunction})). We have%
\begin{equation*}
\nabla w=k\alpha \varphi _{1}{}^{\alpha -1}\nabla \varphi _{1}
\end{equation*}

\begin{equation*}
\Delta w=k\alpha (\alpha -1)\varphi _{1}{}^{\alpha -2}\left\vert \nabla
\varphi _{1}\right\vert ^{2}-\lambda _{1}k\alpha \varphi _{1}{}^{\alpha }.
\end{equation*}%
Thus, by assumption (H$_{2}$), over the ring $\overline{\Omega -K}$ there
exists $\varepsilon >0$ such that 
\begin{equation*}
\varepsilon =\min_{\overline{\Omega -K}}((\alpha -1)\left\vert \nabla
\varphi _{1}\right\vert ^{2}-\lambda _{1}\varphi _{1}^{2}).
\end{equation*}%
Then we have 
\begin{equation*}
-\Delta w-f(x)\leq -\varepsilon k\alpha \varphi _{1}{}^{\alpha -2}+M\varphi
_{1}{}^{\alpha -2}=[M-\varepsilon k\alpha ]\varphi _{1}{}^{\alpha -2}.
\end{equation*}%
On the other hand, the inequality

\begin{equation*}
w(x)\leq C^{+}\text{ on }\partial K
\end{equation*}%
holds if we take, for instance, 
\begin{equation*}
k=\frac{C^{+}}{\max_{x\in \partial K}\{\varphi _{1}(x)^{\alpha }\}},
\end{equation*}%
since $\max_{\partial K}\varphi _{1}{}^{\alpha }>0.$ Then, from the
definition of $M$ given in (\ref{Choice M}) we have that $M-\varepsilon
k\alpha =0$ and thus

\begin{equation*}
\left\{ 
\begin{array}{lr}
-\Delta w(x)-f(x)\leq 0\leq -\Delta u(x)-f(x) & \text{in }\Omega \setminus K,
\\ 
w\leq u & \text{on }\partial K \\ 
w=u=0 & \text{on }\partial \Omega ,%
\end{array}
\right. \text{ }
\end{equation*}%
(the trace, on $\partial K,$ of the very weak solution $u$ is well defined
by the regularity results of \cite{Diaz-Rako}). Then, since $\Delta w\in
L^{1}(\Omega :\delta ),$ we can apply the comparison principle on the ring $%
\Omega \setminus K$ (see, e.g., \cite{Brez-Caz et al}) and we get that 
\begin{equation*}
0<w(x)\leq u(x),\text{ a.e. }x\in \Omega \setminus K,
\end{equation*}%
which ends the proof of (A).

\bigskip

\noindent The proof of the second conclusion will follows again from the
Green's formula. We recall that now $u\in W_{0}^{1,1}(\Omega )$ and $\Delta
u\in L^{1}(\Omega )$, which gives a meaning to the Green's formula since $%
\frac{\partial u}{\partial n}\in L^{1}(\partial \Omega )$. From (\ref{Flat
solu assumpt}), integrating in (\ref{ProbLineal}), we get 
\begin{equation*}
\int_{\partial \Omega }\frac{\partial u}{\partial n}(\sigma )d\sigma
=-\int_{\Omega }\Delta u(x)dx=\int_{\Omega }f(x)dx=0,
\end{equation*}%
and since we already know that $\frac{\partial u}{\partial n}(\sigma )\leq 0$
on $\partial \Omega $ we conclude that $\frac{\partial u}{\partial n}=0$ on $%
\partial \Omega .$

\bigskip \noindent

\noindent \textit{Proof of Lemma} \ref{Lemm uniform y}. The conclusion can
be obtained in some different ways. We will give here some direct arguments.
By applying the Uniform Hopf Inequality to problem (\ref{Auxiliary problems}%
) we get 
\begin{equation*}
\varsigma _{y}(x)\geq C(\Omega )\left( \int_{B_{\varrho }(y)}\delta
(x)dx\right) \delta (x)\text{ a.e. }x\in \Omega .
\end{equation*}%
The function 
\begin{equation*}
\mu (y)=\int_{B_{\varrho }(y)}\delta (x)dx
\end{equation*}%
is a continuous function on the compact $\partial K.$ Thus, there exists $%
\underline{\mu }>0$ such that $\mu (y)\geq \underline{\mu }$ for any $y\in
\partial K$ and then we get 
\begin{equation*}
c_{K}^{\ast }\delta (x)\leq \varsigma _{y}(x)\text{ a.e. }x\in \Omega ,\text{
for any }y\in \partial K,
\end{equation*}%
with 
\begin{equation*}
c_{K}^{\ast }=C(\Omega )\underline{\mu }.
\end{equation*}

\noindent On the other hand, let $\eta (x)$ satisfying 
\begin{equation}
\left\{ 
\begin{array}{ll}
-\Delta \eta =\chi _{\Omega ^{+}} & \text{in }\Omega , \\ 
\eta =0 & \text{on }\partial \Omega .%
\end{array}%
\right.
\end{equation}%
Then, by the comparison principle and well-known regularity results we get

\begin{equation*}
\varsigma _{y}(x)\leq \eta (x)\leq \overline{C}\delta (x)\text{, a.e. }x\in
\Omega ,
\end{equation*}%
for some $\overline{C}>0$ and for any $y\in \partial K$. Thus $\underset{
y\in \partial K}{\sup }C_{K}(y)\leq \overline{C}$.

\bigskip

As in the one-dimensional problem, the existence of nonnegative solutions
with compact support in a larger domain $\widetilde{\Omega }\varsupsetneq
\Omega $, is a simple consequence of the existence of flat solutions on the
small domain.

\begin{corollary}
\label{Corol compact support onedim copy(1)}Let $\widetilde{f}\in L^{1}(%
\widetilde{\Omega })$ be the extension of a given function $f\in
L^{1}(\Omega )$, i.e. such that 
\begin{equation*}
\widetilde{f}(x)=\left\{ 
\begin{array}{ll}
f(x) & \text{if }x\in \Omega , \\ 
0 & \text{if }x\in \widetilde{\Omega }\setminus \Omega .%
\end{array}
\right.
\end{equation*}%
Assume that $f$ satisfies the conditions (H$_{1}$), (H$_{2}$),\textit{\ and
( \ref{Flat solu assumpt}). Let }$u$ be the unique solution of (\ref%
{ProbLineal}) and let $\widetilde{u}$ be the extension of $u$ defined as%
\begin{equation*}
\widetilde{u}(x)=\left\{ 
\begin{array}{ll}
u(x) & \text{if }x\in \Omega , \\ 
0 & \text{if }x\in \widetilde{\Omega }\setminus\Omega .%
\end{array}
\right.
\end{equation*}%
Then $\widetilde{u}$ is the unique weak solution of the problem 
\begin{equation}
\left\{ 
\begin{array}{ll}
-\Delta \widetilde{u}=\widetilde{f}(x) & \text{in }\widetilde{\Omega }, \\ 
\widetilde{u}=0 & \text{on }\partial \widetilde{\Omega }.%
\end{array}
\right.
\end{equation}
\end{corollary}

\begin{remark}
\textrm{The comments on the optimality of the non-negative character of $%
f(x) $ in previous statements of the strong maximum principle made in the
previous section apply also for N-dimensional linear problems. Nevertheless,
the extension to other second order elliptic operators is a delicate point.
The case of Lipschitz coefficients can be treated \ (see Remarks \ref{Rm
First order} and \ref{Rem p-Laplace and obstacle problem}) but, even if $%
f(x) $ is non-negative, there are several counterexamples in the literature
showing that the strong maximum principle may fail for merely bounded
coefficients (see, e.g., \cite{BerdanDRako}, \cite{Apushkinskaya-Nazarov}
and its references). It would be interesting to see if the application of
the Green's function allows to get sharp conditions on the datum $f(x)$ as
in the one-dimensional case (Remark \ref{Rm Rakotoson}).$\hfill%
\mbox{$\quad{}_{\Box}$}$ }
\end{remark}

\begin{remark}
\textrm{\label{Rm Geometrical )} The first condition in (H$_{2}$) has a
geometrical meaning. It requires that if $x_{M}\in \Omega $ is such $\varphi
_{1}(x_{M})=1$ then $x_{M}\in K\subset \Omega ^{+}$. When $\Omega $ is a
ball, $\Omega =B_{R}(0)$, this condition requires that $0\in K\subset \Omega
^{+}$. This holds in the case of Theorem \ref{Thm onedimensional}. Notice
also that the value of $\alpha >1$ which makes possible this condition must
increase if $d(\Omega \setminus\Omega ^{+},0)$ decreases. $\hfill%
\mbox{$\quad{}_{\Box}$}$ }
\end{remark}

\begin{remark}
\textrm{Arguing as in \cite{Alvarez-D}, it seems possible to get a
\textquotedblleft local\textquotedblright\ condition on $f(x)$ on some
neighborhood of some points $x_{0}\in \partial \Omega $ in order to
construct suitable local supersolutions implying that $\frac{\partial u}{%
\partial n}(x_{0})=0$ (see also \cite{Cirmi-Diaz2}). $\hfill%
\mbox{$\quad{}_{\Box}$}$ }
\end{remark}

\begin{remark}
\textrm{\label{Rm Schroedinger}It is quite easy to adapt the above proof to
the case of the Schr\"{o}dinger operator with an absorption potential $V(x)$%
\begin{equation}
\left\{ 
\begin{array}{ll}
-\Delta u+V(x)u\geq f(x) & \text{in }\Omega , \\ 
u=0 & \text{on }\partial \Omega .%
\end{array}
\right.
\end{equation}%
The easy case concerns bounded absorption potentials 
\begin{equation*}
0\leq V(x)\leq M_{V}\text{ on }\Omega .
\end{equation*}%
Indeed, the second part of the proof consists in finding a subsolution on
the region $\Omega -K$ and using again $w(x)=k\varphi _{1}(x)^{\alpha _{V}}$
we see that 
\begin{equation*}
-\Delta w+V(x)w\geq -k\alpha _{V}\varphi _{1}{}^{\alpha _{V}-2}\left(
(\alpha _{V}-1)(\left\vert \nabla \varphi _{1}\right\vert ^{2}-\varphi
_{1}{}^{2}[\lambda _{1}+V(x)]\right) .
\end{equation*}%
It is not difficult to see that we can take $\alpha _{V}\geq \alpha $ (the
exponent when there is no absorption term). Then, $w_{V}(x)=k\varphi
_{1}(x)^{\alpha _{V}}\leq k\varphi _{1}(x)^{\alpha }$, i.e., the subsolution
in the case with absorption is smaller than the one for the case without
absorption. }

\textrm{\noindent The case of an unbounded potential requires additional
approximation arguments. For some related results see \cite{D Ambiguity Free
b}, \cite{D SEMA}, \cite{D Correction SEMA}, \cite{Orsina-Ponce}, \cite%
{DGomezRakoTemam}, \cite{DGomezRako(postTemam)}, \cite{DH Portugalia} and 
\cite{DH linearizacion}. $\hfill\mbox{$\quad{}_{\Box}$}$ }
\end{remark}

\begin{remark}
\textrm{\label{Rm Nonlocal}The case of the Schr\"{o}dinger equation for a
nonlocal diffusion (the fractional Laplace operator) and an unbounded
absorption potential can be also considered with the techniques of the paper 
\cite{DGomezVazquez} (exposition made in \cite{Diaz Schrodingrelativist}). $%
\hfill\mbox{$\quad{}_{\Box}$}$ }
\end{remark}

\begin{remark}
\textrm{\label{Rm First order}The study of the presence of linear transport
terms in the equation can be also carried out (see some similar techniques
in \cite{DGomezRakoTemam}, \cite{DGomezRako(postTemam)} and \cite{BocardoDG}%
). When the negative part of $f(x)$ takes place on an interior subregion of $%
\Omega $ (as indicated in Remark \ref{Rem f negative interior}) and $%
f^{-}(x) $ satisfies some suitable conditions, it is possible to prove the
positivity of the solution of 
\begin{equation}
\left\{ 
\begin{array}{ll}
-\Delta u+\mathbf{b}(x)\cdot \nabla u\geq f(x) & \text{in }\Omega , \\ 
u=0 & \text{on }\partial \Omega ,%
\end{array}
\right.
\end{equation}%
by means of probabilistic methods (see \cite{GDiaz-Diaz Proba}), which also
allow to consider the case of a general elliptic operator with Lipschitz
coefficients. $\hfill\mbox{$\quad{}_{\Box}$}$ }
\end{remark}

\begin{remark}
\textit{\label{Rem p-Laplace and obstacle problem}Theorem \ref{Thm
N-dimensional} can be extended also to the case of some nonlinear diffusion
operators of the type of the p-Laplace operator and the one arising in the
study Bingham fluids and double phase operators (see \cite{Cirmi-Diaz2} and
the presentation made in \cite{D-Pitman} and \cite{Cirmi-D}). Some related
papers applying other properties of the distance function are \cite{Baldelli}
and \cite{Filippuci}. The positivity of solutions is also very relevant in
the class of the obstacle problems. The unilateral condition }$u\geq 0$%
\textit{\ leads to the existence of solutions with compact support once that
the right hand side is assumed to be negative in a neighborhood of the
boundary (see, }$e.g.$\textit{, \cite{brezis ushpeki} and \cite{D-Pitman}).
Theorem \ref{Thm N-dimensional} allows to get the optimality of the
assumptions made on }$f(x)$\textit{\ in order to get solutions with compact
support. Finally, it seems possible to modify the proof of Theorem \ref{Thm
N-dimensional} (specially its first step) as to be applied to suitable fully
nonlinear equations of Monge-Amp\`{e}re type (see \cite{GDiaz-Diaz Monge})
and also to the case of a general elliptic operator with Lipschitz
coefficients. A paper, now in preparation, will develop all these
extensions. }
\end{remark}

\begin{remark}
\textrm{\ \label{Rem Siravov}After the publication of a first version of
this paper in Arxiv, Boyan Sirakov communicated to the authors that the
positivity of solutions and the non-vanishing of the normal derivative could
be proved once the right-hand side has a positive part that is sufficiently
larger than its negative part, in a suitable integral sense. This follows
from Theorem 1.1 of his paper \cite{Sirakov}, dealing with viscosity
solutions of a general class of non-divergent second order elliptic
equations. Such an application is not explicitly mentioned and existence of
flat solutions is not considered in his paper. }
\end{remark}

\begin{remark}
\textrm{\ \label{Rm Romeo} Also, after the publication of a first version of
this paper in Arxiv, Romeo Leylekian communicated to the authors that the
results of this section can be applied to get some new properties of
solutions of higher order elliptic equations involving the bi-Laplacian
operator. }
\end{remark}

\section{Application to some \textit{sublinear} \textit{indefinite}
semilinear equations}

Let us give a short application of the results in previous sections to some 
\textit{sublinear} \textit{indefinite} semilinear equations (see, e.g. \cite%
{Bandle}, \cite{Hernandez-Mancebo-Vega}) and \cite{DiazHern-Ilyas2023}) for
many references and additional results): consider nonnegative \ solutions of
the problem 
\begin{equation}
\left\{ 
\begin{array}{ll}
-\Delta u=\lambda u+m(x)u^{\alpha } & \text{in }\Omega , \\ 
u=0 & \text{on }\partial \Omega ,%
\end{array}%
\right.  \label{Prob1.1}
\end{equation}%
where $\Omega $ is a smooth bounded domain in $\mathbb{R}^{N}$, $0<\alpha <1$%
, $m\in L^{\infty }(\Omega )$ changes sign on $\Omega $ and $\lambda $ is a
real parameter. Problems of this type arise in the study of problems with
nonlinear diffusion in mathematical biology and porous media. \ For
simplicity we will assume 
\begin{equation}
\left\{ 
\begin{array}{c}
m\in L^{\infty }(\Omega )\text{, and }\left\vert \Omega ^{+}\right\vert
,\left\vert \Omega ^{-}\right\vert >0,\text{ with } \\ 
\\ 
\Omega ^{+}=\{x\in \Omega \text{ }\mid m(x)>0\}\text{ and }\Omega
^{-}=\{x\in \Omega \text{ }\mid m(x)<0\}%
\end{array}%
\right.  \label{Hypo Indefinite case}
\end{equation}

We have:

\begin{theorem}
\label{Thm3}Assume (\ref{Hypo Indefinite case}), $m(x)$ satisfies (H$_{1}$),
(H$_{2}$) and let $0\leq \lambda <\lambda _{1}$. Then there is a solution $%
u_{\lambda }>0$ to (\ref{Prob1.1})$.$
\end{theorem}

\noindent \emph{Proof.} We will follow an idea already presented in Lemma
4.2 in \cite{Hernandez-Mancebo-Vega}\textit{. }Let now\textit{\ }$U$ be the
unique solution of the linear problem

\begin{equation}
\left\{ 
\begin{array}{ll}
-\Delta U=m(x) & \text{in }\Omega , \\ 
U=0 & \text{on }\partial \Omega .%
\end{array}
\right.
\end{equation}%
By Theorem \ref{Thm N-dimensional} we know that%
\begin{equation}
U>0\text{ in }\Omega .\text{ }
\end{equation}%
Let us check that function $u_{0}=\left[ (1-\alpha )U\right] ^{1/(1-\alpha
)} $ is a subsolution of problem (\ref{Prob1.1}). Indeed%
\begin{equation*}
\nabla u_{0}=\left[ (1-\alpha )U\right] ^{\alpha /(1-\alpha )}\nabla U
\end{equation*}%
and 
\begin{equation*}
\hbox{div }(\nabla u_{0})=\alpha \left[ (1-\alpha )U\right] ^{(2\alpha
-1)/(1-\alpha )}\left\vert \nabla U\right\vert ^{2}-m(x)\left[ (1-\alpha )U %
\right] ^{\alpha /(1-\alpha )},
\end{equation*}%
which gives%
\begin{eqnarray*}
-\Delta u_{0}-\lambda u_{0}-m(x)(u_{0})^{\alpha } &=&-\alpha \left[
(1-\alpha )U\right] ^{\frac{(2\alpha -1)}{(1-\alpha )}}\left\vert \nabla
U\right\vert ^{2} \\
-\lambda \left[ (1-\alpha )U\right] ^{1/(1-\alpha )} &\leq &0
\end{eqnarray*}%
if $\lambda \geq 0.$

\noindent On the other hand, as a supersolution $u^{0}=C\psi $ we pick the
function $\psi >0$ as the solution of the problem%
\begin{equation}
\left\{ 
\begin{array}{ll}
-\Delta \psi =\lambda \psi +1 & \text{in }\Omega , \\ 
\psi =0 & \text{on }\partial \Omega ,%
\end{array}%
\right.
\end{equation}%
(recall that $0\leq \lambda <\lambda _{1}$). Then we have 
\begin{equation*}
-\Delta u^{0}-\lambda u^{0}-m(x)(u^{0})^{\alpha }=C-C^{\alpha }m(x)\psi
^{\alpha }=C^{\alpha }(C^{1-\alpha }-m(x)\psi ^{\alpha })>0
\end{equation*}%
if $C>(\left\Vert m\right\Vert _{L^{\infty }}\left\Vert \psi \right\Vert
_{L^{\infty }}^{\alpha })^{1/(1-\alpha )}.$ On the other hand, we know that $%
U,\psi \in C_{0}^{1}(\overline{\Omega })$, with 
\begin{equation*}
\frac{\partial U}{\partial n}\leq 0\text{ and }\frac{\partial \psi }{%
\partial n}<0\text{ on }\partial \Omega .
\end{equation*}%
Then, for $C$ large enough we know that $u_{0}\leq u^{0}$ on $\Omega $. Then
the existence of a solution $u_{\lambda }>0$ to (\ref{Prob1.1}) follows from
the super and subsolution method in \cite{Hernandez-Mancebo-Vega}. $\hfill %
\mbox{$\quad{}_{\Box}$}$

\begin{remark}
\textrm{Many other results on the problem (\ref{Prob1.1}) can be found in 
\cite{Hernandez-Mancebo-Vega}. For instance, if in Theorem \ref{Thm3} we
assume, in addition, that the unique weak solution $U$ of the linear problem
(\ref{ProbLineal}) satisfies $\frac{\partial U}{\partial n}<0$ on $\partial
\Omega $ then there is uniqueness of positive solutions to problem (\ref%
{Prob1.1}). Indeed, the uniqueness of positive solutions can be obtained by
means of a suitable change of unknowns (see Theorem 4.4 of \cite%
{Hernandez-Mancebo-Vega} and the extension presented in \cite%
{DiazHern-Ilyas2023}), or by means of some hidden convexity arguments (see,
e.g., \cite{Diaz PAFA}).$\hfill\mbox{$\quad{}_{\Box}$}$ }
\end{remark}

\begin{remark}
\textrm{It is easy to see that the above subsolution becomes, in fact, an
exact solution, if we replace the coefficient $m(x)$ by the function 
\begin{equation*}
\widehat{m}(x)=m(x)+\alpha (1-\alpha )^{\frac{(2\alpha -1)}{(1-\alpha )}} 
\frac{\left\vert \nabla U\right\vert ^{2}}{U}+\lambda (1-\alpha
)^{1/(1-\alpha )}U.
\end{equation*}%
In that case, when $\alpha \in (\frac{1}{2},1)$ the positive solution
becomes a \textquotedblleft flat solution\textquotedblright\ since in that
case the exponent of $U$ in the expression of $\nabla u_{0}$ is $\alpha
/(1-\alpha )>1.\hfill\mbox{$\quad{}_{\Box}$} $ }
\end{remark}

\section{An application to the linear parabolic problem}

Another application of Theorems \ref{Thm onedimensional} and \ref{Thm
N-dimensional} deals with the question of the positivity of solutions of the
linear parabolic problem%
\begin{equation}
\left\{ 
\begin{array}{ll}
u_{t}-\Delta u=f(x,t) & \text{in }\Omega \times (0,+\infty ), \\ 
u=0 & \text{on }\partial \Omega \times (0,+\infty ), \\ 
u(x,0)=u_{0}(x) & \text{on }\Omega .%
\end{array}%
\right.  \label{Parabolic problem}
\end{equation}%
By using the arguments for stationary equations of Sections 2 and 3, jointly
to some related results (D\'{\i}az-Fleckinger \cite{DF}), we can prove that,
under suitably changing sign conditions on $u_{0}(x)$ and/or on $f(x,t)$
(even if it occurs for any $t>0$), we get the global positivity of $u(x,t)$
on $\Omega $, for large values of time. Here is a special statement of this
kind of results:

\bigskip

\begin{theorem}
Let $\Omega $ as Theorem \ref{Thm N-dimensional} and let $f\in
L_{loc}^{1}(0,+\infty :L^{1}(\Omega :\delta ))$ with 
\begin{equation*}
f(x,t)\geq g(x)\text{ a.e. }x\in \Omega ,\text{ a.e. }t>0
\end{equation*}%
with $g(x)$ changing sign and satisfying the conditions of Theorem \ref{Thm
N-dimensional}. Assume $u_{0}\in L^{1}(\Omega :\delta )$ such that there
exists a $\widehat{u_{0}}\in L^{2}(\Omega )$, with $u_{0}(x)\geq $ $\widehat{%
u_{0}}(x)$ a.e. $x\in \Omega $ such that 
\begin{equation*}
\int_{\Omega }\widehat{u_{0}}\varphi _{1}=0.
\end{equation*}%
Let $u\in C([0,+\infty ):L^{1}(\Omega :\delta ))$ be the unique mild
solution of problem (\ref{Parabolic problem}). Then, there exists $t_{0}\geq
0$ such that $u(x,t)>0$ a.e. $x\in \Omega ,$ for any $t>t_{0}$.
\end{theorem}

\noindent \emph{Proof.} Let $v(x)$ be the solution of the stationary problem 
\begin{equation}
\left\{ 
\begin{array}{ll}
-\Delta v=g(x) & \text{in }\Omega , \\ 
v=0 & \text{on }\partial \Omega ,%
\end{array}%
\right.
\end{equation}%
and let $w\in C([0,+\infty ):L^{2}(\Omega ))$ be the unique solution of the
homogeneous parabolic problem%
\begin{equation}
\left\{ 
\begin{array}{ll}
w_{t}-\Delta w=0 & \text{in }\Omega \times (0,+\infty ), \\ 
w=0 & \text{on }\partial \Omega \times (0,+\infty ), \\ 
w(x,0)=\widehat{u_{0}}(x) & \text{on }\Omega .%
\end{array}%
\right.
\end{equation}%
Then, it is clear that the function 
\begin{equation*}
\underline{u}(t,x)=v(x)+w(t,x)
\end{equation*}%
is a subsolution to problem (\ref{Parabolic problem}). Moreover, by
Proposition 2.3 of \cite{DF} (originally, an unpublished result due to the
first author and J. M. Morel) 
\begin{equation}
\left\vert w(t,x)\right\vert \leq C\left\Vert \widehat{u_{0}}\right\Vert
_{L^{2}(\Omega )}e^{-\lambda _{2}t}\varphi _{1}(x)\text{, a.e. }x\in \Omega ,%
\text{ for any }t>t_{0},  \label{Estimate heat equation}
\end{equation}%
where $\lambda _{2}$ is the second eigenvalue of the Laplacian operator on $%
\Omega .$ Then, the conclusion follows by Theorem \ref{Thm N-dimensional}
and (\ref{Estimate heat equation}) since we have the comparison principle in
the class of mild solutions $u\in C([0,+\infty ):L^{1}(\Omega :\delta ))$ on
the Banach space $X=L^{1}(\Omega :\delta )$ (see, e.g. \cite{Brez-Caz et al}%
).$\hfill \mbox{$\quad{}_{\Box}$}$

\section*{Acknowledgments}

The authors would like to thank Jean Michel Rakotoson, Boyan Sirakov and
Romeo Leylekian for some comments on a previous version of this paper (see,
Remarks \ref{Rm Rakotoson}, \ref{Rem Siravov} and \ref{Rm Romeo},
respectively). 

\bigskip 

All authors declare that they have no conflicts of interest.

\noindent \vskip1cm 
\flushleft{\small
\begin{tabular}{ll}\small
J.I. D\'{\i}az & J. Hernández \\
Instituto Matemático Interdisciplinar (IMI) & Instituto Matemático Interdisciplinar (IMI) \\
Dept. Análisis Matemático y Matem\'atica Aplicada & \\
U. Complutense de Madrid & U. Complutense de Madrid \\
Parque de las Ciencias & Parque de las Ciencias\\
28040-Madrid. Spain & 28040-Madrid. Spain\\
{\tt jidiaz@ucm.es} & {\tt jesus.hernande@telefonica.net}
\end{tabular}
 }

\end{document}